\documentclass[reqno,10pt]{amsart}
\usepackage{amscd,amssymb,amsmath,color,mathtools}
\usepackage{paralist}
\usepackage{latexsym}
\usepackage[all]{xy}
\usepackage[colorlinks,allcolors=blue]{hyperref}
\usepackage[T1]{fontenc}
\usepackage{color}
\usepackage[english]{babel}
\usepackage{textcomp}
\usepackage{amsmath}
\usepackage{amssymb}
\usepackage{graphicx}
\usepackage{algorithm}
\usepackage{algpseudocode}

\newcommand{\lyxmathsym}[1]{\ifmmode\begingroup\def\b@ld{bold}
  \text{\ifx\math@version\b@ld\bfseries\fi#1}\endgroup\else#1\fi}

\def\..{{,\dots,}}
\def\:{{\colon}}

\def\int{{\rm int}}

\theoremstyle{definition}

\usepackage{booktabs} % Top and bottom rules for tables

\graphicspath{{figures/}} % Location of the graphics files
 \usepackage{tikz} 
 \usetikzlibrary{calc}

\usepackage{tikz}
\usetikzlibrary{arrows}

\begin{document}

\author{Or Raz}

\date{\today}

\title{Independent sets of non-geometric lattices and the maximal adjoint}

\address{Einstein Institute of Mathematics, The Hebrew University of Jerusalem, Giv'at Ram, Jerusalem, 91904, Israel}
\email{or.raz1@mail.huji.ac.il}

\begin{abstract}
We construct a family of independent sets for finite, atomic, and graded lattices, extending the well-known cryptomorphism between geometric lattices and matroids. This construction leads to an embedding theorem into geometric lattices that preserves the set of atoms. We then apply these results to adjoint matroids, providing new characterizations of adjoints and partially proving a conjecture on the combinatorial derived matroid. Finally, we use our characterization of adjoints to compute the adjoint lists of several simple examples.
\end{abstract}

\maketitle

\section{Introduction}
The construction of adjoint matroids was introduced in 1974 by ALAN L. C. CHEUNG in \cite{[2]} as an attempt to extend the duality principle for projective space, i.e., the duality between points and hyperplanes. In the same vein, the set of matroids with an adjoint can be viewed as a generalization of the set of linear matroids. In recent years there has been renewed interest in adjoints coming from linear codes and dependencies between circuits of a matroid. As shown in \cite{[2]}, \cite{[3]}, \cite{[4]} and \cite{[5]} there does not always exist an adjoint, and even if it does exist it may not be unique.

The case of representable matroids (over an arbitrary field) was developed simultaneously in two separate works \cite{[7]} and \cite{[8]}, showing (in part) that every representable matroid has a representable adjoint corresponding to its representation. In \cite{[6]} the combinatorial derived matroid was introduced in an attempt to produce a generalization of an adjoint to arbitrary matroids. The first characterization of an adjoint of arbitrary matroids using hyperplanes was given in \cite{[12]}, which was used in the classification of adjoint sequences.

The foundation of our new approach to adjoints was initially developed in \cite{[9]}, where a family of independent sets of finite, atomic, and graded lattices was introduced and leveraged to prove an embedding theorem for lattices in which the family of independent sets constitutes a matroid. In this paper, we extend these results to cases where the family of independent sets does not constitute a matroid, allowing for further applications to the problem of adjoints.

The structure of this paper is as follows. Section 2 is dedicated to preliminaries on lattices ,matroids, and the adjoint matroid. In Section 3 we discuss the family of independent sets of a finite, atomic and graded lattice, recalling results from \cite{[9]} and presenting a new embedding theorem. In Section 4 we show the effects of basic lattice operations on its family of independent sets and a class of special maximal independent sets. Section 5 moves to working with adjoint matroids centered around the characterization theorems of adjoints as an application of the family of independent sets. We further demonstrate how these methods can be applied to establish known results, specifically proving that every representable matroid and co-rank $3$ matroid possess an adjoint in a significantly simpler manner. Additionally, an immediate consequence is the identification of the explicit unique adjoint of a modular matroid. In Section 6 we present some partial results showing that the combinatorial derived matroid is an adjoint and discussing the existence of a maximal adjoint. Finally, this paper also includes an appendix in which we formulate an algorithm to check for adjoints based on the theorems in Section 5. We then use this algorithm to calculate the list of adjoints of some matroids of small rank and co-rank. \\

Throughout the paper, we use lower case letters when referring to elements, upper case letters when referring to sets, upper case letters in calligraphic font when referring to families of sets, and upper case letters in gothic font when referring to collections of families. For example, we may have $x\in X\in \mathcal{X}\in \mathfrak{X}$.

\section{Preliminaries}
We start with some basic lattice properties; for a more extensive background on matroids and geometric lattices, see \cite{[1]}. A partially ordered set $\mathcal{L}$ in which every pair of elements has a unique supremum (also called join and denoted $\vee$) and a unique infimum (also called meet and denoted $\wedge$) is called a lattice. Two important elements, which always exist in a finite lattice, are the greatest element (denoted $\hat{1}$) and the least element (denoted $\hat{0}$). Given two elements $A\leq B$ in a lattice $\mathcal{L}$, we define the segment $[A,B]_\mathcal{L}=\left\{C\mid A\leq C\leq B\right\}$. Notably, a segment of a lattice forms a lattice itself.

\subsection{Definition:}
\label{Definition1}
The dual lattice $\mathcal{L}^{opp}$ of $\mathcal{L}$ is defined as the lattice with the same underlying set and the reverse order. \\

We say $B$ covers $A$ in $\mathcal{L}$, and write $A\lessdot B$ if $A<B$ and for all $C$ such that $A\leq C\leq B$ we have $C=B$ or $C=A$. $A$ is an atom if it covers $\hat{0}$ and a coatom if it is covered by $\hat{1}$. A finite lattice is called atomic if every element is a join of atoms, and coatomic if its dual lattice is atomic. We denote the set of atoms of $\mathcal{L}$ by $A(\mathcal{L})$.

Finally, we say that a finite lattice is graded if there exists a rank function $\rho:\mathcal{L}\rightarrow \mathbb{N}$ with the following properties:

\begin{enumerate}
\item $\rho(\hat{0})=0$.
\item Compatibility with the order of $\mathcal{L}$, $A\leq B\iff \rho(A)\leq \rho(B)$.
\item Compatibility with the covering relation of $\mathcal{L}$, $A\lessdot B\Rightarrow \rho(A)= \rho(B)+1$.
\end{enumerate}

\subsection{Definition:}
\label{Definition2}
A geometric lattice is a finite graded lattice equipped with a submodular rank function, meaning it satisfies the inequality $\rho(A)+\rho(B)\leq \rho(A\vee B)+\rho(A\land B)$. \\

Next, we turn to some fundamental properties of matroids; for a more detailed discussion, refer to \cite{[1]}.

We adopt standard conventions: a matroid $\mathcal{M}$ refers to its set of bases. We denote the family of independent sets by $I(\mathcal{M})$, the collection of circuits by 
 $\mathcal{C}(\mathcal{M})$, the associated geometric lattice by $\mathcal{L}(\mathcal{M})$, and the ground set by $E(\mathcal{M})$. When the reference to a specific matroid is clear, we omit $\mathcal{M}$, using, for instance, $\mathcal{L}$ instead of $\mathcal{L}(\mathcal{M})$.\\

 An element $a\in E(\mathcal{M})$ will be called a loop of $\mathcal{M}$ if $\left\{a\right\}\notin I(\mathcal{M})$. Two elements $a,b\in E$ will be called parallel if both are not loops and $\left\{a,b\right\}\notin I(\mathcal{M})$ . We will say that a matroid is simple if it has no loops or parallel elements.
 
\subsection{Definition:}
\label{Definition3}
Given two matroids $\mathcal{M}_1, \mathcal{M}_2$ on the ground set $E$, we say $\mathcal{M}_1\leq \mathcal{M}_2$ in the weak order of matroids if $I(\mathcal{M}_1)\subseteq I(\mathcal{M}_2)$.\\

The dual matroid of $\mathcal{M}$, denoted $\mathcal{M}^*$ is the matroid defined on the same ground set, where its bases (maximal independent sets) are precisely the complements of the bases of $\mathcal{M}$.
The geometric lattice associated with $\mathcal{M}^*$ is denoted by $\mathcal{L}^*$ and a loop of $\mathcal{M}^*$ is referred to as a co-loop.\\

A known property is that the complements of the hyperplanes (another name for co-atoms of geometric lattices) are precisely the circuits of the dual matroid. This property enables us to interpret an extension of the dual lattice of a geometric lattice, as discussed in the following sections, as a matroid with ground set $\mathcal{C}(\mathcal{M})$. We structure our definitions accordingly to reflect this perspective:

\subsection{Definition:}
\label{Definition4}
Let $\mathcal{L}, \mathcal{L}^{\triangle}$ be geometric lattices. We say $\mathcal{L}^{\triangle}$ is an adjoint of $\mathcal{L}^*$ if there exists a rank preserving order embedding $f: \mathcal{L}^{opp} \rightarrow \mathcal{L}^{\triangle}$ which restricts to a bijection from $ A(\mathcal{L}^{opp})$ to $ A(\mathcal{L}^{\triangle})$. For matroids $\mathcal{M},\mathcal{M}^\triangle$, we say that $\mathcal{M}^\triangle$ is an adjoint of $\mathcal{M}$ if $\mathcal{L}(\mathcal{M}^\triangle)$ is an adjoint of $\mathcal{L}(\mathcal{M}^*)$. \\

The following two properties are known for adjoints (for example, in \cite{[5]}):

\begin{enumerate}
\item $A\lessdot B\in \mathcal{L}\Rightarrow f(A)\lessdot f(B)$
\item $f(A\land B)=f(A)\land f(B)$

\end{enumerate}

\section{The family of independent sets and embedding theorems}

We start with the definition of $NBB$ sets as given in \cite{[10]}.
\subsection{Definition:}
\label{Definition5}
Let $(\mathcal{L},\omega)$ be a pair of a finite lattice $\mathcal{L}$ and a partial order $\omega$ on its set of atoms $A(\mathcal{L})$. A nonempty set $D \subseteq A(\mathcal{L})$ of atoms is bounded below or $BB$ if there exist $a\in A(\mathcal{L})$ such that $a$ is a strict lower bound for all $d\in D$ in the order $\omega$ and $a\leq \vee D$ in $\mathcal{L}$.\\

A set $B \subseteq A(\mathcal{L})$ is called $NBB$ if $B$ does not contain any bounded below subset. 

Given a finite lattice $\mathcal{L}$ and $P$ the set of all linear orders on $A(\mathcal{L})$ we define the family of independent sets $I(\mathcal{L})$ of $\mathcal{L}$ as follows:

\[
I\left(\mathcal{L}\right)=\left\{ A\in 2^{A(\mathcal{L})}\mid\exists\omega\in P\text{ s.t. }A\text{ is an }NBB\text{ set in }\left(\mathcal{L},\omega\right)\right\}
\]

We refer to the elements in $I(\mathcal{L})$ as independent sets of $\mathcal{L}$. Before continuing, we note the following immediate properties:
\begin{enumerate}
    \item A subset of an independent set is itself independent. In other words, independent sets form an abstract simplicial complex.
    \item A set containing a single atom is not $BB$ in any order.
    \item sets containing at most two atoms are always independent. This is done by choosing any order where one of the two elements is first.
\end{enumerate} 

\subsection*{Example:} Let $\mathcal{L}\subseteq 2^{[4]}$ be the non-geometric lattice containing all subsets of size $0,1,4$ and $\left\{1,2\right\}$, ordered by inclusion. Following Definition \ref{Definition5}, we see that $I(\mathcal{L})$ contains $\left\{\left\{1,2,3\right\},\left\{1,2,4\right\}\right\}$ and all subsets of size $\leq 2$.\\

The constructions in this section can be thought of as a generalization of the following observation, proven in \cite{[9]}.

\subsection{Theorem:}
\label{Theorem1}
If $\mathcal{L}$ is a geometric lattice, then $I(\mathcal{L})$ is the family of independent sets of the matroid corresponding to $\mathcal{L}$.\\

Next, we review key properties and definitions established in \cite{[9]}.

\subsection{Definition:}
\label{Definition6}
 Let $\mathcal{L}$ be a finite, atomic, and graded lattice, and $A\subseteq A(\mathcal{L})$. If $\vee A=X\in \mathcal{L}$ we say that $A$ spans $X$ and define $\rho(A)\coloneqq \rho(X)$. An independent set $A$ of $\mathcal{L}$ is called geometric if $\rho(A)=\left|A\right|$.\\

\subsection{Lemma:}
\label{Lemma1}
Let $\mathcal{L}$ be a finite, and graded lattice (not necessarily atomic), and $A\in I(\mathcal{L})$. There exists an atom $a\in A$ such that $\rho(A)\lneq \rho(A\setminus \left\{a\right\})$.

\subsection{Theorem:}
\label{Theorem2}
Let $\mathcal{L}$ be a finite, atomic, and graded lattice, and $A\in I(\mathcal{L})$. It is always true that $\left|A\right|\leq \rho(A)$. Moreover, if $\left|A\right|\neq \rho(A)$ then there exists an independent set $A\subsetneq B\in I(\mathcal{L})$ such that $\vee B=\vee A$ and $\left|B\right|=\rho(B)$. Consequently, the rank of an element equals the size of a maximal independent set contained in it.\\

The following embedding theorem was introduced in \cite{[9]} and covers the case in which $I(\mathcal{L})$ is a family of independent sets of a matroid.

\subsection{Theorem:}
\label{Theorem3}
Let $\mathcal{L}$ be a finite, atomic and graded lattice. If $I(\mathcal{L})$ is a family of independent sets of a matroid, then there exists a rank-preserving order embedding $f: \mathcal{L}\rightarrow \mathcal{P}$ for which the restriction to the set of atoms $f\mid_{A(\mathcal{L})}:A(\mathcal{L})\rightarrow A(\mathcal{P})$ is a bijection. With $\mathcal{P}$ being the geometric lattice corresponding to $I(\mathcal{L})$.\\

We proceed with the formulation of an embedding theorem for cases where $I(\mathcal{L})$ does not constitute a family of independent sets of a matroid.

\subsection{Lemma:}
\label{Lemma2}
Let $\mathcal{L}$ be a finite and graded lattice (not necessarily atomic) and let $f':A(\mathcal{L}) \rightarrow A(\mathcal{P})$ be a bijection between the sets of atoms of $\mathcal{L}$ and a geometric lattice $\mathcal{P}$. If $f'$ can be extended to a rank-preserving order embedding $f: \mathcal{L}\rightarrow \mathcal{P}$, then $f(I(\mathcal{L}))\subseteq I(\mathcal{P})$. 

\subsection*{Proof:}
Let $f: \mathcal{L}\rightarrow \mathcal{P}$ be such an embedding and assume $f(I(\mathcal{L}))\nsubseteq I(\mathcal{P})$. We can choose a minimal $A\in f(I(\mathcal{L}))\setminus I(\mathcal{P})$. Let $a\in f^{-1}(A)$ with $\rho(f^{-1}(I)\setminus \left\{ a\right\})\lneq \rho(f^{-1}(A))$ guaranteed by Lemma \ref{Lemma1}. As $A\setminus \left\{ f(a)\right\}$ is an independent set of $\mathcal{P}$ and $A$ is not an independent set of $\mathcal{P}$, we must have
\[ f(a)\in \vee_\mathcal{P}(A\setminus \left\{ f(a)\right\})\leq _\mathcal{P} f(\vee_\mathcal{L}(f^{-1}(A)\setminus \left\{ a\right\})) \],
which is a contradiction to the choice of $a$.

\subsection{Theorem:}
\label{Theorem4}
Let $\mathcal{L}$ be a finite, atomic and graded lattice, and let $f':A(\mathcal{L}) \rightarrow A(\mathcal{P})$ be a bijection between the set of atoms of $\mathcal{L}$ and the set of atoms of a geometric lattice $\mathcal{P}$ of the same rank. $f'$ can be extended to a rank-preserving order embedding $f: \mathcal{L}\rightarrow \mathcal{P}$ if and only if $f(I(\mathcal{L}))\subseteq I(\mathcal{P})$ and $f(A)\notin I(\mathcal{P})$ for every $A\in 2^{A(\mathcal{L})}$ such that $\left|A\right|\lneq \rho\left(A\right)$. 

\subsection*{Proof:}
If there exists a function $f$ then Lemma \ref{Lemma2} implies the first condition. The second condition holds as if $f(A)\in I(\mathcal{P})$ then $\rho(f(A))=\left|A\right|$ and so $f$ is rank-preserving.\\
In the other direction, we define $f(X)\coloneqq \vee_\mathcal{P}(f(A(X)))$. First, notice that $\rho(X)\leq \rho(f(X))$ for all $X\in \mathcal{L}$ as $f(X)$ contains the image of a geometric spanning set of $X$. To see $\rho (X)\geq \rho(f(X))$ for all $X\in \mathcal{L}$ observe that it cannot contain a larger independent set of $\mathcal{P}$ as it will contradict the second condition.
  To see $f$ is an embedding, let $X\in \mathcal{L}$, $a\in A(\mathcal{L})\setminus A(X)$ and $B$ be a geometric independent set spanning $X$. We have that $f(A)\cup \left\{ f(a)\right\}\in I(\mathcal{P})$ is an independent set of size $\rho(X)+1$ in $f(X)$ which contradicts $f$ being rank-preserving.

\section{Lattice operations}
In this section we present the effects of basic lattice operations on its family of independent sets. 

Let $\mathcal{L}$ be a finite, atomic, and graded lattice, and let $A\in I(\mathcal{L})$ be an independent set. By definition, $A$ is $NBB$ with respect to some linear order. Notably, we can always choose a linear order that begins with the atoms in $A$. We denote such an order, where $A$ remains $NBB$, by $\omega _A$. Now, let $B\subsetneq A$ such that $A$ is $NBB$ with respect to a linear order $\omega _A$ that starts with the atoms in $B$. We denote this order by $\omega_{A,B}$.

\subsection{Definition:}
\label{Definition7}
Let $\mathcal{L}$ be a finite, atomic and graded lattice. we define the following basic lattice operations:
\begin{enumerate}
\item[Restriction:] Let $X\in \mathcal{L}$ be an element of rank $\geq 1$. The restriction of $\mathcal{L}$ to $X$ is the finite, atomic, and graded lattice $[\hat{0},X]$. For clarity, we define $A(X)\coloneqq A([\hat{0},X])$.
\item[Contraction:] Let $X\in \mathcal{L}$ be an element of rank $\leq \rho(\mathcal{L})-1$. The contraction of $\mathcal{L}$ to $X$ is the finite and graded lattice $[X,\hat{1}]$.
\item[Truncation:]  Let $m\leq \rho(\mathcal{L})$. The truncation of $\mathcal{L}$ by $m$ is the finite, atomic, and graded lattice $T(\mathcal{L},m)=\left\{ X\in \mathcal{L}\mid \rho(X)< m\right\}\cup \left\{\hat{1}_{\mathcal{L}}\right\}$. To see the resulting poset is indeed a lattice observe that a lower set of a finite lattice is a finite meet semilattice, so that adding a top element makes it a lattice.
\item[Dual:  ] The dual $\mathcal{L}^{opp}$ of $\mathcal{L}$ is the finite and graded lattice with the same underlying set and the reverse order.
\end{enumerate}

\subsection{Definition:}
\label{Definition8}
Let $\mathcal{L}$ be a finite, atomic and graded lattice. We say $\mathcal{L}'$ is a minor of $\mathcal{L}$ if it is obtained from $\mathcal{L}$ by a sequence of contractions and restriction.

\subsection{Theorem:}
\label{Theorem5}
The basic lattice operations described above influence the family of independent sets of  $\mathcal{L}$ in the following manner:
\begin{enumerate}
\item $I([\hat{0},X])=\left\{ I\cap A(X)\mid I\in I\left(\mathcal{L}\right)\right\}$.
\item  
$
I([X,\hat{1}])\cong (\left\{ I\mid \exists \omega_{I\cup I_X,I}, \vee I_X=X\right\} \cap\left(\cup_{X\lessdot Y}A\left(Y\right)\right))/_{\sim_{X}}
$, with $\sim_{X}$ being the equivalence relation on $A(\mathcal{L})$ having $a\sim_{X} b$ if there exists $X\lessdot Y$ such that $a,b\in Y\setminus X$. The first set can be best understood by selecting an independent set $I_x$ that spans $X$, then examining all $NBB$ that include $I_X$ while following a linear order that begins with the atoms outside $I_X$. Next, we remove atoms that do not belong to $[X,\hat{1}]$ and quotient by atoms that are smaller than the same covering element.
\item  $I(T(\mathcal{L},m))=$\\ $\left\{ I\in I(\mathcal{L})\mid \rho (I)\leq m \right\}\cup \left\{ I\in I(\mathcal{L})\mid \exists a\in I \text{ such that } \rho (I\setminus \left\{a\right\})\lneq m\right\}$.
\item If $\mathcal{L}$ is coatomic, then a set of coatoms $\mathcal{X}\subset A(\mathcal{L}^{opp})$ is a geometric independent set if and only if there exist a linear order on $\mathcal{X}=\left\{ X_1,...,X_n\right\}$ and a geometric independent set $\left\{ a_1,...,a_{d-1}\right\}$ of $\mathcal{L}$ with $\omega_{\left\{ a_1,...,a_{n}\right\}}$ such that $\land_{i=1} ^t X_i= \vee_{i=1} ^{d-t} a_i$ for all $t\in [n]$.
\end{enumerate}

\subsection*{Proof:}
$(1)$ Given a linear order $\omega$ on $A(\mathcal{L})$, we can derive a linear order on the restriction by considering the induced order $\omega\mid_X$ on $A(X)$. An $NBB$ set $I$ in $\omega$ then induces an $NBB$ set $I\cap A(X)$ in $\omega\mid_X$. Conversely, a linear order on $A(X)$ can be extended to a linear order on $\mathcal{L}$ by appending the remaining atoms in any order at the end. In this induced order, every $NBB$ set of $[\hat{0},X]$ remains an $NBB$ set of $\mathcal{L}$ 

$(2)$ We begin by describing the isomorphism, which maps the atoms of $[X,\hat{1}]$ to their corresponding equivalence classes under $\sim_X$, assigning each $Y$  to the class consisting of atoms in $A(Y)$. Let $I\subseteq A([X,\hat{1}])$ be an $NBB$ set with respect to $\omega$. By selecting an independent set $I_X$ that spans $X$ and arranging it in an $NBB$ order, we ensure that $I\cup I_X$ remains $NBB$ with respect to $\omega_{I\cup I_X,I}$, as required by definition. Conversely, the construction proceeds by eliminating atoms that do not cover $X$ and taking the quotient in $\omega_{I\cup I_X,I}$, yielding $I$, which retains the $NBB$ property within the constructed order.

$(3)$ Given a linear order $\omega$ on $A(\mathcal{L})=A(T(\mathcal{L},m))$, every $NBB$ set in $\mathcal{L}$ with rank $\leq m$ remains an $NBB$ set in $T(\mathcal{L},m)$. If $\rho(I)\gneq m$ and there exists $a\in I$ such that $\rho (I\setminus \left\{a\right\})\lneq m$, then there exists a linear order $\omega_{\left\{a\right\}}$ for which $I$ is $NBB$ in $\mathcal{L}$, and $I$ remains an $NBB$ set in $T(\mathcal{L},m)$ with respect to $\omega_{\left\{a\right\}}$. For the remaining case, where $\rho\gneq m$ and no $a\in I$ satisfies $\rho (I\setminus \left\{a\right\})\lneq m$, any linear order on $A(\mathcal{L})$ can be taken. Letting $a$ be the first atom, we have $a\leq_{\mathcal{L}} \vee (I\setminus \left\{a\right\})$. Thus, $I$ is not $NBB$ with respect to any order in the truncation.

$(4)$ holds as a result of Lemma \ref{Lemma1}. Consider such an order on $\mathcal{X}$  and select any geometric independent set $\left\{ a_1,...,a_{n}\right\}$ with $\omega_{\left\{ a_1,...,a_{n}\right\}}$ that spans $\land_{X\in \mathcal{X}} X$. Since $\mathcal{X}$ is geometric, we have $\rho(\land_{i=1} ^{t+1} X_i)+1=\rho(\land_{i=1} ^t X_i)$ for every $t\in [k-1]$. Thus, we can find atoms $a_{n+1},...a_{d-1}$ extending $\left\{ a_1,...,a_{n}\right\}$ to a spanning set of $\land_{i=1} ^{t-1} X_i$ for each $t$. The other direction follows from:
\[
 n=\rho(\vee_{i=1} ^{n} a_i)=\rho(\land_{i=1} ^n X_i)=\left|\mathcal{X}\right| 
\]

\subsection*{Remark:} We compare the effect of contraction for geometric lattices (or simple matroids) with general finite, graded and atomic lattices. Contracting an element of a non-geometric lattice depends on the ordering of the independent sets, even if the contraction is atomic. One consequence of this difference is that the family of finite, atomic and graded lattices with $I(\mathcal{L})$ a matroid is not closed under the operation of taking minors. For example, the lattice $\mathcal{L}=\mathcal{L}(\mathcal{U}_{6,4})\setminus \left\{\left\{123 \right\},\left\{124 \right\},\left\{126 \right\},\left\{134 \right\},\left\{136 \right\} \right\}$, with $\mathcal{U}_{6,4}$ being the uniform matroid of rank $4$ on the ground set $[6]$. One can easily check that $\mathcal{L}$ is an atomic and graded lattice with $I(\mathcal{L})=I(\mathcal{U}_{6,4})$ (obviously a matroid). However, the contraction $[1,[6]]_\mathcal{L}$, which is still a finite, atomic and graded lattice, has a non-matroid family of independent sets $I([1,[6]]_\mathcal{L})=I(\mathcal{U}_{5,3})\setminus \left\{\left\{234 \right\}, \left\{236 \right\} \right\}$.\\

For the second part of this section we introduce the following special kind of maximal independent set.

\subsection{Definition:}
\label{Definition9}
Let $\mathcal{L}$ be a finite, atomic and graded lattice, a maximal independent set $I\in I(\mathcal{L)}$ will be called a basis of $\mathcal{L}$ if every subset of $I$ is geometric.\\

Our reason for calling these maximal independent sets bases is that, in the same way as bases of matroids, they admit the following basis-exchange property. 

\subsection{Theorem:}
\label{Theorem6}
Let $\mathcal{L}$ be a finite, atomic and graded lattice. If $B$ is a basis of $\mathcal{L}$ and $I$ is a maximal independent set, then for every atom $a\in B$ there exists an atom $b\in I$ such that $(B\setminus \left\{a \right\})\cup \left\{b \right\}$ is a maximal independent set.\\

This property comes from the fact that $\vee(B\setminus \left\{a \right\})\neq \hat{1}$ for every $a\in B$. Another matroid-like property bases possess is the existence of a fundamental circuit.

\subsection{Lemma:}
\label{Lemma3}
Let $B$ be a basis for $\mathcal{L}$. For every $a\notin B$, $B\cup \left\{a \right\}$ contains a unique minimal dependent set denoted $\mathcal{C}_{B,a}$ . 

\subsection*{Proof:}
Assume $C_1,C_2$ are minimal dependent sets contained in $B\cup \left\{a \right\}$ and let $x\in C_1\setminus C_2$. $C_1\setminus \left\{x \right\}$ and $C_1\setminus \left\{a \right\}$ are independent, adding an atom not in the span of an independent set to it is again independent. Therefore, $\vee C_1=\vee(C_1\setminus \left\{x \right\})=\vee(C_1\setminus \left\{a \right\})$ and therefore a geometric independent set. We then have $(B\cup \left\{a \right\})\setminus \left\{x \right\}$ independent since adding each atom in $B\setminus \left\{x \right\}$ to $C_1\setminus \left\{x \right\}$ (in any order) will increase its rank, and as $\vee(C_1\setminus \left\{x \right\})=\vee(C_1\setminus \left\{a \right\})$ this is also true for $C_1\setminus \left\{a \right\}$. We now have $C_2\subseteq (B\cup \left\{a \right\})\setminus \left\{x \right\}$, which contradicts $C_2$ being dependent.\\

We continue our exploration of bases by examining how they transform under the operation of taking duals. Our goal is to define the dual "image" of bases, we observe that each subset $I'$ of size $\rho(\mathcal{L})-1$ of a basis $\mathcal{I}$ corresponds to the unique coatom $\vee I'$. We define $I^{opp}=\left\{\vee (I\setminus \left\{a \right\})\mid a\in I \right\}$ as the basis of $\mathcal{L}^{opp}$ corresponding to $I$. To see $\mathcal{I}^{opp}$ is a basis, observe that for every subset $\mathcal{J}\subsetneq \mathcal{I}^{opp}$ and every $I^{opp}\in \mathcal{I}^{opp}\setminus	\mathcal{J}$ we have $\rho (\vee_{\mathcal{L}^{opp}}\mathcal{J})\lneq \rho(\vee_{\mathcal{L}^{opp}}(\mathcal{J}\cup \left\{I^{opp} \right\}))$.

\subsection{Theorem:}
\label{Theorem7}
Let $\mathcal{B}$ and $\mathcal{B}^{opp}$ be the family of bases of $\mathcal{L}$ and $\mathcal{L}^{opp}$ respectively. If $\mathcal{L}$ is coatomic, then the function $opp:\mathcal{B}\rightarrow \mathcal{B}^{opp}$ with $opp(I)=I^{opp}$ is a bijection.

\subsection*{Proof:}
Let $I,J$ be bases of $\mathcal{L}$ with $opp(I)=opp(J)$. For all $a\in I$ there exists some $b\in J$ with $\vee (I\setminus \left\{a \right\})=\vee (J\setminus \left\{b \right\})$. $opp$ being injective now follows by induction on $I\setminus \left\{a \right\}$ and $J\setminus \left\{b \right\}$ as bases of $[\hat{0},\vee (I\setminus \left\{a \right\})]$. \\
Let $\mathcal{J}$ be a basis of $\mathcal{L}^{opp}$, for every $J\in \mathcal{J}$ we have a unique atom $a=\vee  _{\mathcal{L}^{opp}}(\mathcal{J}\setminus \left\{J \right\})\in A(\mathcal{L})$. Therefore, $a\nleq J$ and $a\leq I$ for every $I\in \mathcal{J}\setminus	\left\{J \right\}$. Consequently, $I=\left\{\vee  _{\mathcal{L}^{opp}}(\mathcal{J}\setminus \left\{J \right\})\mid J\in \mathcal{J} \right\}$ is a basis of $\mathcal{L}$ such that $opp(I)=J$.\\

\section{Adjoint matroids}
In this section, we extend the construction developed in previous sections to the study of adjoint matroids. We begin by recalling a well-established result that links the coatoms—also referred to as hyperplanes—of a geometric lattice with the circuits of its dual, as presented in \cite{[1]}.

\subsection{Lemma:}
\label{Lemma4}
Let $\mathcal{M}$ be a matroid and $\mathcal{M}^*$ its dual matroid. The circuits of $\mathcal{M}^*$ are exactly the complement of the coatoms of $\mathcal{L}(\mathcal{M})$.

We continue with a definition of an important type of dependent sets related to $(\mathcal{L}^*)^{opp}\coloneqq(\mathcal{L}(\mathcal{M})^*)^{opp}$, the dual lattice of the dual matroid. We have already shown in Lemma \ref{Lemma1} that the size of an independent set of a lattice is $\leq$ its rank, therefore, every set larger than its rank must be dependent. We denote the complement of this family of sets by:

\[
\mathcal{S}(\mathcal{M})=\left\{ \mathcal{S}\in 2^{\mathcal{C}(\mathcal{M})}\mid\left|\mathcal{S}\right|\leq \rho_{(\mathcal{L}^*)^{opp}}\left(\cap_{C\in\mathcal{S}}C^c\right)\right\}
\]

This definition is purposely written in a complicated way, to show that it depends only on the matroid $\mathcal{M}$. Thinking of circuits of $(\mathcal{M})$ as atoms of $(\mathcal{L}^*)^{opp}$, and forgetting its complicated connection to $\mathcal{M}$, we are describing sets of atoms not larger than their rank.  

\subsection*{Example:} Let $\mathcal{U}_{3,3}$ be the uniform matroid of rank
$3$ on the ground set $[6]$. Its circuits are all the subsets of $[6]$ of size $4$ and therefore the atoms of $(\mathcal{L}(\mathcal{U}_{3,3})^*)^{opp}$ are the sets of size $2$. We already know that pairs of atoms are always independent, so every pair is in $\mathcal{S}(\mathcal{U}_{3,3})$. We have $\rho ((\mathcal{L}(\mathcal{U}_{3,3})^*)^{opp})=3$, and so sets of four atoms are never in $\mathcal{S}({\mathcal{U}_{3,3}})$. Sets of three atoms $\left\{a,b,c\right\}$ are in $\mathcal{S}({\mathcal{U}_{3,3}})$ if and only if $\left|a\cup b\cup c\right|\geq 4$.\\

%The second type of dependent sets we define are called "evenly dependent", denoted by:

%\[
%\mathcal{S}_m(\mathcal{M})= \left\{ rank_{\mathcal{M}}\left(\cap_{C\in\mathcal{S}}C^{c}\right)=rank_{\mathcal{M}}\left(\cap_{C\in\mathcal{S}\setminus\left\{ A\right\} }C^{c}\right),\,\forall A\in\mathcal{S}\right\} 
%\]

%Notice that every circuit of $(\mathcal{L}^*)^{opp}$ is evenly dependent.

It was proven (for example, in "On Adjoints and Dual Matroids") that not all matroids have an adjoint and even if they have one it may not be minimal. We can reformulate Theorem \ref{Theorem4} to give a necessary and sufficient condition for a matroid to be an adjoint of another matroid:

\subsection{Theorem:}
\label{Theorem8}
Let $\mathcal{M}$ be a matroid of rank $d$ on the ground set $[n]$ and $\mathcal{M}^{\triangle}$ a matroid of rank $n-d$ on the ground set $\mathcal{C}(\mathcal{M})$. $\mathcal{M}^{\triangle}$  is an adjoint of $\mathcal{M}$ if and only if $I((\mathcal{L}^*)^{opp})\subseteq I(\mathcal{M}^{\triangle})\subseteq \mathcal{S}(\mathcal{M)}$.

\begin{enumerate}
\item $\mathcal{M}^{\triangle}$  is an adjoint of $\mathcal{M}$.
\item  $I((\mathcal{L}^*)^{opp})\subseteq I(\mathcal{M}^{\triangle})\subseteq \mathcal{S}(\mathcal{M)}$.
\item $\left\{ C\mid i\notin C\in\mathcal{C}\left(\mathcal{M}\right)\right\}$ is a hyperplane of $\mathcal{M}^{\triangle}$ for every $i\in[n]$ not a co-loop.  
\end{enumerate}

\subsection*{Proof:}

\begin{enumerate}
\item[$(1\Rightarrow 2)$] If $\mathcal{M}^{\triangle}$ is an adjoint of $\mathcal{M}$ then by Lemma \ref{Lemma2} we have $I((\mathcal{L}^*)^{opp})\subseteq \mathcal{M}^{\triangle}$. Let $\mathcal{A}\notin \mathcal{S}(\mathcal{M})$ and $I$ be a geometric independent set of $(\mathcal{L}^*)^{opp}$ spanning $\cap_{C\in \mathcal{A}}C^c$. By Theorem \ref{Theorem4}, we have $\left|I\right|=\rho_{\mathcal{M}^\triangle}(\mathcal{A})< \left|\mathcal{A}\right|$, which means $\mathcal{A}\notin I(\mathcal{M}^\triangle)$.
\item[$(2\Rightarrow 3)$] If $i\in E$ is not a co-loop, then there exists a basis $i\notin B\in \mathcal{M}$. Examining the set of fundamental circuits $C_B\setminus \left\{C_{B,i}\right\}$ of $B$, we obtain an independent set of rank $n-d-1$ contained in $\left\{ C\mid i\notin C\in\mathcal{C}\left(\mathcal{M}\right)\right\}$. It is not of full rank as, by the definition of $I((\mathcal{L}^*)^{opp})$, it will give us an independent set of rank $n-d+1$. By the same reasoning, we must have $\left\{ C\mid i\notin C\in\mathcal{C}\left(\mathcal{M}\right)\right\}$ closed.
\item [$(3\Rightarrow 1)$] Since geometric lattices are co-atomic, the set $\left\{ C\mid i\notin C\in\mathcal{C}\left(\mathcal{M}\right)\right\}$ corresponds to the hyperplanes of $(\mathcal{L}^*)^{opp}$. Consequently, $(\mathcal{L}^*)^{opp}$ must be order-embedded into $L(\mathcal{M}^{\triangle})$. Furthermore, as $A(L(\mathcal{M}^{\triangle}))$ forms a partition of $\mathcal{C}(\mathcal{M})$ and $\mathcal{M}^{\triangle}$ has rank $n-d$ , the embedding preserves order and becomes a bijection when restricted to $A((\mathcal{L}^*)^{opp})$.
\end{enumerate}

\subsection{Corollary:}
It follows that if $I((\mathcal{L}^*)^{opp})$ or $\mathcal{S}(\mathcal{M)}$ are matroids, then $\mathcal{M}$ has an adjoint. Moreover, $I((\mathcal{L}^*)^{opp})$ serves as the minimal adjoint, while $\mathcal{S}(\mathcal{M)}$ represents the maximal adjoint within the weak order on the class of matroids with ground set $\mathcal{C}(\mathcal{M})$.\\

Two separate works, \cite{[8]} by James Oxley and Suijie Wang and \cite{[7]} by Relinde Jurrius and Ruud Pellikaan, define an adjoint for linear matroids. We will use the notation of Oxley and Wang as Jurrius and Pellikaan used a dual definition. 

Let $\mathcal{M}$ a representable matroid on the ground set $E=(e_1,...,e_n)$, and let $\varphi :E\rightarrow \mathbb{F}^m$ be a representation of $\mathcal{M}$. To each circuit $C$ of $\mathcal{M}$ is associated with respect to $\varphi$ a unique vector (up to scalar multiplication) $v_C=(c_1,...,c_n)\in \mathbb{F}^m$ such that $\sum_{i=1}^{n}c_{i}\varphi\left(e_{i}\right)=0$, where $c_i\neq 0$ if and only if $i\in C$.

\subsection{Definition:}
\label{Definition10}
To a pair of representable matroids and representation $(\mathcal{M},\varphi)$ define the derived matroid $(\delta \mathcal{M}, \delta \varphi)$ with ground set $\mathcal{C}\left(\mathcal{M}\right)$ to be the representable matroid determined by $\delta \varphi (C)=v_C$.\\

In the work of Jurrius and Pellikaan it was proven that derived matroids are adjoints; here we give an alternative proof.  

We first show that a collection of circuits $\mathcal{A}=\left\{C_1,...,C_k\right\}$ of a linear matroid $\mathcal{M}$ corresponding to an independent set of $I((\mathcal{L}^*)^{opp})$ is independent in the derived matroid, for any representation. As $\mathcal{A}$ corresponds to an independent set, we can linearly order its circuit such that:
\[
C_i\setminus (\cup_{j<i}C_j)\neq \emptyset, \forall i\in [k]
\]
We then have $\mathcal{A}$ an independent set in the derived matroid by simple linear algebra. \\

To see sets not in $\mathcal{S}(\mathcal{M})$ are dependent in $\delta \mathcal{M}$, we use the following observations:

\begin{enumerate}
    \item If $\mathcal{A}$ is a minimal set of circuits not in $\mathcal{S}(\mathcal{M})$ then $\left|\cup_{C\in \mathcal{A}}C\right|= \rho_\mathcal{M}(\cup_{C\in \mathcal{A}}C)+\left|\mathcal{A}\right|+1$. Also, for any $C\in \mathcal{A}$ we have $\left|\cup_{C\in \mathcal{A}}C\right|=\left|\cup_{C\in \mathcal{A}\setminus \left\{C'\right\}}C\right|$.
     \item If $C,C'$ are two circuits such that $\left|C\cup C'\right|=\rho_\mathcal{M}(C\cup C')+1$ then for every $e\in C\cap C'$ there exists a unique circuit $C''$ such that $C''\subseteq C\cup C'\setminus \left\{e\right\}$ and $C,C',C''$ are dependent in $\delta \mathcal{M}$.
\end{enumerate}

We proceed by induction on the number of circuits in $\mathcal{A}$ , with the base case established by the second observation. Let $\mathcal{A}$ be independent and minimal set of circuits not contained in $\mathcal{S}(\mathcal{M})$. Using the second observation, we eliminate an element $e\in C\in \mathcal{A}$ by replacing each occurrence of $e\in C'\neq C$ with $C''$. This process results in an independent set $\mathcal{A}'$ of identical size and rank. Since $e\in C\setminus (\cup _{C'\in \mathcal{A}'\setminus \left\{C\right\}}C)$, we observe that $\mathcal{A}'\setminus \left\{C\right\}$ still contains a minimal set of circuits not in $\mathcal{S}(\mathcal{M})$. The proof then follows by the induction hypothesis.\\

An established result we can revisit is that every co-rank $3$ matroid $\mathcal{M}$ has an adjoint. To see this, consider $\mathcal{M}'= \left\{ \left\{C_1,C_2,C_3\right\}\mid \cup_{i=1}^3C_i=[n]\right\}$. A straightforward verification shows that $\mathcal{M}'=\mathcal{S}(\mathcal{M})$, confirming that it serves as an adjoint.\\

Our next objective is to simplify the computation of adjoint requirements, given that evaluating $I((\mathcal{L}^*)^{opp})$ can be computationally demanding. To begin, we observe that the fundamental circuit of any basis in $(\mathcal{L}^*)^{opp}$ must be a circuit in every adjoint of $\mathcal{M}$, as it is not contained in $\mathcal{S}(\mathcal{M})$. Additionally, since $\mathcal{L}^*$ is a geometric lattice, every basis of $(\mathcal{L}^*)^{opp}$ takes the form $(E\setminus	B)^{opp}=\left\{C_{B,e}\mid e\notin B\right\}=\mathcal{C}_B$ for a basis $B$ of $\mathcal{M}$. The following provides an explicit representation of the "fundamental circuits" of $(\mathcal{L}^*)^{opp}$:

\subsection{Lemma:}
\label{Lemma5}
Let $B$ be a basis of $\mathcal{L}$ and $C\in \mathcal{C}(\mathcal{M})$, then:
\[
\mathcal{C}_{\mathcal{C}_{\mathcal{B}},C}=\left\{C_{B,e}\mid e\in C\setminus B\right\}\cup \left\{C\right\}
\]
\subsection*{Proof:}
If $e\in C\setminus B$ then as $e\notin C_{B,a}$ for every $a\neq e$ we must have $(\mathcal{C}\setminus (\left\{C_{B,e}\right\}))\cup \left\{C\right\}$ independent in $(\mathcal{L}^*)^{opp}$ and so $C_{B,e}\in \mathcal{C}_{\mathcal{C}_{\mathcal{B}},C}$. On the other hand, let $e\notin C\setminus B$ and assume by negation $C_{B,e}\in \mathcal{C}_{\mathcal{C}_{\mathcal{B}},C}$. We then have $\mathcal{C}_{\mathcal{C}_{\mathcal{B}},C}\setminus \left\{C_{B,e}\right\}$ independent and $e$ not contained in any of its elements, consequently $\mathcal{C}_{\mathcal{C}_{\mathcal{B}},C}$ is independent, contradicting it being a circuit.\\

We now present an algorithm for computing maximal independent sets of $(\mathcal{L}^*)^{opp}$ based on its bases. Beginning with a basis $\mathcal{C}_B$ of $(\mathcal{L}^*)^{opp}$, we have $\rho(\mathcal{C}_B\setminus \left\{C_{B,i_1}\right\})=n-d-1$ for every $i_1\notin B$. Consequently, $(\mathcal{C}_B\setminus \left\{C_{B,i_1}\right\})\cup \left\{C_1\right\}\in I((\mathcal{L}^*)^{opp})$ for any circuit $C_1$ of $\mathcal{M}$ that contains $i_1$. This process can be iterated: in the $k$'th step, we replace $C_{B,i_k}$ with $C_k$, where $i_k\in C_k\setminus (B\cup (\cup_{j=1}^{k-1}C_j))$. As a result, we obtain $(\mathcal{C}_B\setminus \left\{C_{B,i_1},...,C_{B,i_k}\right\})\cup \left\{C_1,..,C_k\right\}\in I((\mathcal{L}^*)^{opp})$.

\subsection{Lemma:}
\label{Lemma6}
Every maximal independent set of $I((\mathcal{L}^*)^{opp})$ is obtained by the above algorithm.

\subsection*{Proof:}
Let $\mathcal{I}$ be a maximal independent set, and consider an ordered sequence $\left\{C_1,...,C_{n-d}\right\}$ of its elements for which $\omega_\mathcal{I}$ exist. By Theorem \ref{Theorem7}, there is a basis $B$ of $\mathcal{M}$ such that $E\setminus	B=\left\{a_1,...,a_{d-n}\right\}$, satisfying $\cap_{i=0} ^k C_i^c= \vee_{i=1} ^{n-d-k} a_i$ for all $t\in \left\{0,..,d-n\right\}$, where the join of the $a_i$'s is taken in $\mathcal{L}^*$. Additionally, we define $\cap_{i=0} ^0C_i=E$. It follows that $\left\{a_{d-n-k-1},...,a_{d-n}\right\}\in C_k$ and $\left\{a_1,...,a_{k}\right\}\notin C_k$ for all $k\in [d-n-1]$. Consequently, the algorithm initiates with $\mathcal{C}_B$ and systematically replaces circuits in the given order $(a_{d-n},...a_1)$ with the circuits $\left\{C_1,...,C_{n-d}\right\}$, ultimately constructing $\mathcal{I}$.

We can now formulate the following equivalent condition to being an adjoint:

\subsection{Theorem:}
\label{Theorem9}
$\mathcal{M}^\triangle$ is an adjoint of $\mathcal{M}$ if and only if $\mathcal{C}_{B}$ forms a basis of $\mathcal{M}^\triangle$ for every basis $B$ of $\mathcal{M}$, and $\mathcal{C}_{\mathcal{C}_{B},C}$ forms a circuit of $\mathcal{M}^\triangle$ for every $C\notin \mathcal{C}_{B}$.

\subsection*{Proof:}
We have already established that $\mathcal{C}_{B}$ forms a basis of $\mathcal{M}^\triangle$. To see $\mathcal{C}_{\mathcal{C}_{B},C}$ forms a circuit of $\mathcal{M}^\triangle$, note that $\mathcal{C}_{\mathcal{C}_{B},C}\notin \mathcal{S}(\mathcal{M})$ and every proper subset of $\mathcal{C}_{\mathcal{C}_{B},C}$ remains independent in $(\mathcal{L}^*)^{opp}$.\\
Conversely, since $\mathcal{C}_{B}$ is a basis of $\mathcal{M}^\triangle$, we ensure that $\mathcal{M}^\triangle$ has the required rank. Next, we verify that condition $(3)$ of Theorem \ref{Theorem8} holds. Let $e$ be a non-coloop element of $\mathcal{M}$ and $B$ a basis of $\mathcal{M}$ that does not contain $e$. We aim to show that $\vee\left\{C\mid e\notin C(\mathcal{M})\right\}=\vee(\mathcal{C}_{B}\setminus \left\{C_{B,e}\right\})$ which then forms a hyperplane of $\mathcal{M}^\triangle$. If $e\in C \in \mathcal{C}(\mathcal{M})$, then $(\mathcal{C}_{B}\setminus \left\{C_{B,e}\right\})\cup \left\{C\right\}\in I((\mathcal{L}^*)^{opp})$ and so $C\notin \vee(\mathcal{C}_{B}\setminus \left\{C_{B,e}\right\})$. If $e\notin C\in \mathcal{C}(\mathcal{M})$, then by Lemma \ref{Lemma5} $C_{B,e}\notin \mathcal{C}_{\mathcal{C}_{B},C}$ and therefore $C\in \vee(\mathcal{C}_{B}\setminus \left\{C_{B,e}\right\})$.\\

\subsection{Corollary:}
Every family of circuits of $\mathcal{M}$ not contained in $\mathcal{S}(\mathcal{M})$ must contain a fundamental circuit of $(\mathcal{L}^*)^{opp}$.\\
 \\
 
\subsection{Theorem:}
\label{Theorem10}
If $\mathcal{M}$ is a matroid corresponding to a modular geometric lattice $\mathcal{L}$, then $I((\mathcal{L}^*)^{opp})$ is its only adjoint matroid.

\subsection*{Proof:}
As $\mathcal{L}$ is modular we have every maximal independent set of $I((\mathcal{L}^*)^{opp})$ of the form $\mathcal{C}_{B}$ for a basis $B$ of $\mathcal{M}$. If $\mathcal{M}^\triangle$ is an adjoint of $\mathcal{M}$ with a basis $\mathcal{B}\notin I((\mathcal{L}^*)^{opp})$, we must have some basis $B$ of $\mathcal{M}$ with $C_1\in \mathcal{C}_{B}$ such that for every $C_2\in \mathcal{B}$, $(\mathcal{C}_{B}\setminus \left\{C_{1}\right\})\cup \left\{C_{2}\right\}\notin I((\mathcal{L}^*)^{opp})$. Therefore, $(\mathcal{C}_{B}\setminus \left\{C_{1}\right\})\cup \left\{C_{2}\right\}\in (\mathcal{S}(\mathcal{M}))^c$ contradicting $\mathcal{M}^\triangle$ being an adjoint.

\section{The Combinatorial derived matroid and maximal adjoints}

Applying the characterizations of adjoints developed in section 5 we can show that the combinatorial derived matroid defined in \cite{[6]} is an adjoint for two relatively simple cases. We repeat some of the definitions from \cite{[6]} for the sake of completeness.

\subsection{Definition:}
\label{Definition11}
Let $\mathcal{M}$ be a matroid, $\mathfrak{A}\subset 2^{C(\mathcal{M})}$ a collection of circuit sets, and $k\in \mathbb{N}$. We define the following operations:

\begin{enumerate}
\item 
\[
\mathfrak{A}_{k}=\left\{ \mathcal{A}\in\mathfrak{A}\mid\left|\mathcal{A}\right|=k\right\}, \mathfrak{A}_{\leq k}=\left\{ \mathcal{A}\in\mathfrak{A}\mid\left|\mathcal{A}\right|\leq k\right\}  
\]
\item 
\[
\epsilon\left(\mathcal{\mathfrak{A}}\right)=\mathcal{\mathfrak{A}}\cup\left\{ \left(\mathcal{A}_{1}\cup\mathcal{A}_{2}\right)\setminus\left\{ C\right\} \mid\mathcal{A}_{1},\mathcal{A}_{2}\subseteq\mathcal{\mathfrak{A}},\,\mathcal{A}_{1}\cap\mathcal{A}_{2}\notin\mathfrak{A},\,C\in\mathcal{A}_{1}\cap\mathcal{A}_{2}\right\} 
\]
\item 
\[
\mathfrak{\uparrow\mathcal{\mathfrak{A}}}=\left\{ \mathcal{A}\in P\left(C\left(\mathcal{M}\right)\right)\mid\exists\mathcal{A}'\in\mathfrak{A},\,\mathcal{A}'\subseteq\mathcal{A}\right\} 
\] 
\end{enumerate}

In \cite{[6]} the second and third operations were used to construct a matroid on the set of circuits $C(\mathcal{M})$. Given some initial collection $\mathfrak{  A}\subseteq 2^{C(\mathcal{M})}$ Define $\mathfrak{A}_0=\mathfrak{A}$ and construct inductively the collections $\mathfrak{A}_{i+1}=\uparrow \epsilon \mathfrak{A}_{i}$. As this sequence is increasing and contained in the finite set $2^{C(\mathcal{M})}$, we have the well-defined limit:

\[
D(\mathfrak{A})=\cup_{i\geq 0} \mathfrak{A}_i
\]

The next theorem explains the abuse of notation in writing $D(\mathfrak{A})$, showing that it is indeed the family of dependent sets of the unique matroid constructed from $\mathfrak{A}$.

\subsection{Theorem} (Proposition 4.7. of \cite{[6]})
\label{Theorem11}
If $\emptyset\notin \mathfrak{A}$ then $D(\mathfrak{A})$ is the family of dependent sets of a matroid on $C(\mathcal{M})$.

\subsection{Definition:}
\label{Definition12}
The combinatorial derived matroid $\delta \mathcal{M}$ of $\mathcal{M}$ is defined as the matroid on the set of circuits $C(\mathcal{M})$ with dependent sets $D(2^{C(\mathcal{M})}\setminus \mathcal{S}(\mathcal{M}))$.\\

In \cite{[6]} they did not have the criteria for $\delta \mathcal{M}$ being an adjoint of $\mathcal{M}$. Using Theorem \ref{Theorem8} we know it is enough for $I((\mathcal{L}^*)^{opp})$ to be independent in $\delta \mathcal{M}$ as $2^{C(\mathcal{M})}\setminus \mathcal{S}(\mathcal{M})$ is already dependent. In particular, if $D(2^{C(\mathcal{M})}\setminus \mathcal{S}(\mathcal{M}))= \\  \uparrow (2^{C(\mathcal{M})}\setminus \mathcal{S}(\mathcal{M}))$ then $\delta \mathcal{M}$ is an adjoint. Two simple examples follow: 

\subsection{Theorem}
\label{Theorem12}
If $\mathcal{M}$ is of co-rank $3$ then $\delta \mathcal{M}$ is an adjoint of $\mathcal{M}$.

\subsection*{Proof:}
The only elements in $2^{C(\mathcal{M})}\setminus \mathcal{S}(\mathcal{M})$ of size three are sets of the form $\left\{ C_{1},C_{2},C_{3}\right\}$ with:
\[
1=\rho_{\mathcal{M}^*}(C_{1}^c\cap C_{2}^c\cap C_{3}^c)=\rho_{\mathcal{M}^*}(C_{1}^c\cap C_{2}^c)=\rho_{\mathcal{M}^*}(C_{1}^c\cap C_{3}^c)=\rho_{\mathcal{M}^*}(C_{2}^c\cap C_{3}^c)
\]

A set of size three in $\epsilon (2^{C(\mathcal{M})}\setminus \mathcal{S}(\mathcal{M}))$ is of the form $\left\{C_{2},C_{3},C_{4}\right\}$ with some $C_1\in \mathcal{C}(\mathcal{M})$ such that $\left\{ C_{1},C_{2},C_{3}\right\},\left\{ C_{1},C_{2},C_{4}\right\} \in 2^{C(\mathcal{M})}\setminus \mathcal{S}(\mathcal{M})$, which is again in $2^{C(\mathcal{M})}\setminus \mathcal{S}(\mathcal{M})$. As we also have $(\uparrow (2^{C(\mathcal{M})}\setminus \mathcal{S}(\mathcal{M})))_{3}=(2^{C(\mathcal{M})}\setminus \mathcal{S}(\mathcal{M}))_{3}$ and the maximal $\rho_{\mathcal{M}^*}$ of a set is $3$ we get $D(\delta \mathcal{M})=2^{C(\mathcal{M})}\setminus \mathcal{S}(\mathcal{M})$,

\subsection{Theorem} 
\label{Theorem13}
$\delta U(k,n)$ is an adjoint of $U(k,n)$, the uniform matroid of rank $k$ on the ground set $[n]$.

\subsection*{Proof:}
Let $\mathcal{A}_1,\mathcal{A}_2\in 2^{C(U(k,n))}\setminus \mathcal{S}(U(k,n))$ with $C\in \mathcal{A}_1\cap \mathcal{A}_2$. As $U(k,n)^*=U(n-k,n)$ we have:
\[
n-k-\rho_{\mathcal{M}^*}(\cap_{A\in\mathcal{A}_i}A)=\left|\left(\cup_{A\in\mathcal{A}_i}A\right)\setminus C\right|+1<\left|\mathcal{A}_{i}\right|
\]
Therefore, we must also have:
\[
n-k-\rho_{\mathcal{M}^{*}}(\cap_{A\in\left(\mathcal{A}_{1}\cup\mathcal{A}_{2}\right)\setminus\left\{ C\right\} }A)\leq n-k-\rho_{\mathcal{M}^{*}}(\cap_{A\in\mathcal{A}_{1}\cup\mathcal{A}_{2}}A)
\]
\[
=\left|\left(\cup_{A\in\mathcal{A}_{1}\cup\mathcal{A}_{2}}A\right)\setminus C\right|+1\leq\left|\left(\cup_{A\in\mathcal{A}_{1}}A\right)\setminus C\right|+1+\left|\left(\cup_{A\in\mathcal{A}_{2}}A\right)\setminus C\right|+1-1
\]
\[
<\left|\mathcal{A}_{1}\right|+\left|\mathcal{A}_{2}\right|-1=\left|\left(\mathcal{A}_{1}\cup\mathcal{A}_{2}\right)\setminus\left\{ C\right\} \right|
\]

To complete the proof, we observe that if $\mathcal{A}_1,\mathcal{A}_2\in \uparrow (2^{C(U(k,n))}\setminus \mathcal{S}(\mathcal{M}))$ then $\left(\mathcal{A}_{1}\cup\mathcal{A}_{2}\right)\setminus\left\{ C\right\} $ will contain either a subset in $2^{C(U(k,n))}\setminus \mathcal{S}(U(k,n))$ already contained in one of the $\mathcal{A}_{i}$'s or contain a collection of the form $\left(\mathcal{B}_{1}\cup\mathcal{B}_{2}\right)\setminus\left\{ C\right\} $ with $\mathcal{B}_1,\mathcal{B}_2\in 2^{C(U(k,n))}\setminus \mathcal{S}(U(k,n))$, which is again in $2^{C(U(k,n))}\setminus \mathcal{S}(U(k,n))$. Consequently, 
 $D(\delta U(k,n))=\uparrow (2^{C(U(k,n))}\setminus \mathcal{S}(U(k,n)))$.\\
%Using the above construction and theorem we are now ready to construct the maximal adjoint of a matroid.

%\subsection{Definition:}
%Let $\mathcal{M}$ be a matroid of rank $d$ on the ground set $[n]$. The maximal adjoint $\mu \mathcal{M}$ of $\mathcal{M}$ is a matroid on the ground set $C(\mathcal{M})$ and the following collection of dependent sets:

%\[
%D(\mu \mathcal{M})=\ \uparrow D^{(n-d)}(2^{C(\mathcal{A})}\setminus \mathcal{S}(\mathcal{M}))
%\]
%with $D^{(0)}(\mathfrak{A})=\mathfrak{A}$ and $D^{(k)}(\mathfrak{A})=(D(D^{(k-1)}(\mathfrak{A})\cup \mathfrak{A}))_{\leq k}$ for all $k\in [n-d]$.\\

%Those of you who have read \cite{[13]} will notice the similarities between the maximal adjoint and the derive combinatorial matroid. The only difference being the construction of dependent sets of increasing size separately with the purpose of not combining dependent sets with an intersection that will become dependent at a later iteration. It may be the case that these definitions are equivalent, or that the derived combinatorial matroid is a non-maximal adjoint. This leads us to the main theorem of this section:

We notice that in all observed cases, as suggested in \cite{[6]}, we have $\delta \mathcal{M}$ a maximal adjoint of $\mathcal{M}$ in the weak order of matroids. Using the tools developed in \cite{[11]}, we can construct a maximal matroid structure on a set of circuits if some conditions hold. We start by introducing the construction, let $\mathcal{X}\subseteq P(E)$ for the ground set $E$. A $\mathcal{X}$-matroid is a matroid $\mathcal{M}$ on the ground set $E$ such that $\mathcal{X}\subseteq C(\mathcal{M})$. \cite{[11]} introduced the following upper bound on the rank function on a $\mathcal{X}$ -matroid. 

\subsection{Definition:}
\label{Definition13}
A sequence $\mathcal{A}=(X_1,...,X_k)$ will be called a proper $\mathcal{X}$-sequence if $X_i\in \mathcal{X}$ for all $i\in [k]$ and $X_i\subsetneq \cup_{j=1}^{i-1}X_j$ for all $2\leq i\leq k$. For $F\subseteq E$, we define $val(F,A)=\left|F\cup(\cup_{i=1}^{k}X_{i})\right|-k$.

\subsection{Lemma:}
\label{Lemma7}
Let $\mathcal{M}$ be an $\mathcal{X}$-matroid and $F\subseteq E$. Then $\rho_\mathcal{M}(F)\leq val(F,S)$ for any proper $\mathcal{X}$-sequence $\mathcal{A}$.  Furthermore, if the equality holds, then $\rho_\mathcal{M}(F\setminus \left\{e\right\})=\rho_\mathcal{M}(F)-1$ for all $e\in F\setminus \cup_{X\in \mathcal{A}}X$ and $\rho_\mathcal{M}(F\cup \left\{e\right\})=\rho_\mathcal{M}(F)$ for all $e\in \cup_{X\in \mathcal{A}}X$.\\

Using this upper bound, the following function was introduced and shown to be a candidate for a maximal matroid on the set of circuits.

\subsection{Theorem:}
\label{Theorem14}
Let $\mathcal{X}$ be a family of subsets of a ground set $E$ and the set of $\mathcal{X}$-matroids not empty. If the following function $val_{\mathcal{X}}:2^{E}\rightarrow\mathbb{Z}$ is submodular, then it is the rank function of the unique maximal $\mathcal{X}$-matroid on the ground set $E$.

\[
val_{\mathcal{X}}\left(F\right)=min\left\{ val\left(F,\mathcal{A}\right)\mid\mathcal{A}\text{ is a proper }\mathcal{X}\text{-sequence}\right\} 
\]
\\
 In our case, we saw in Theorem \ref{Theorem9} that every collection of circuits of the form $\mathcal{C}_{\mathcal{C}_B,C}$ is a circuit for every adjoint matroid, therefore the natural choice is $\mathcal{X}=\cup_{B\in\mathcal{M}} \left\{ \mathcal{C}_{\mathcal{C}_{\mathcal{B}},C}\mid C\notin\mathcal{C}_{\mathcal{B}}\right\}$ making every adjoint of $\mathcal{M}$ an $\mathcal{X}$-matroid. Moreover, if $B\in \mathcal{M}$ and $\mathcal{A}=\left\{C\notin \mathcal{C}_B\right\}$, we have:
 \[
val_{\mathcal{X}}\left(C(\mathcal{M})\right)\leq val(C(\mathcal{M}),\mathcal{A})=\left|C(\mathcal{M})\right|-\left|C(\mathcal{M})\setminus \mathcal{C}_B\right|=\left|\mathcal{C}_B\right|=\rho(\mathcal{M}^*)
\]
 Therefore $val_{\mathcal{X}}\left(\mathcal{A}\right)\leq \rho(\mathcal{M}^*)$ for every $\mathcal{A}\subseteq C(\mathcal{M})$.\\

\subsection{Lemma:}
\label{Lemma8}
$val_{\mathcal{X}}\left(\mathcal{D}\right)<\left|\mathcal{D}\right|$ for every $\mathcal{D}\in \uparrow (2^{C(\mathcal{M})}\setminus \mathcal{S}(\mathcal{M}))$.
\subsection*{Proof:}
As $val_{\mathcal{X}}(\mathcal{A})\leq val_{\mathcal{X}}(\mathcal{A}\cup\left\{C\right\})\leq val_{\mathcal{X}}(\mathcal{A})+1$ for every $\mathcal{A}\in 2^{C(\mathcal{M})}$ and $C\in C(\mathcal{M})$ it is enough to prove the inequality for $\mathcal{D}\in (2^{C(\mathcal{M})}\setminus \mathcal{S}(\mathcal{M}))$. If $\rho_{\mathcal{M}^*}(\cap_{C\in \mathcal{D}}C^c)=0$, we must have $\left|\mathcal{D}\right|>\rho(\mathcal{M}^*)$ and the lemma follows from the above discussion. Otherwise, there exists an independent set $I$ of $\mathcal{M}^*$ of size $\rho_{\mathcal{M}^*}(\cap_{C\in \mathcal{D}}C^c)$ such that $I\cap C=\emptyset$ for all $C\in\mathcal{D}$. Let $B$ be a basis of $\mathcal{M}$ contained in the complement of $I$ and $\mathcal{A}=\left\{ \mathcal{C}_{\mathcal{C}_{B},C}\mid C\in\mathcal{D}\right\}$. As $\mathcal{A}$ is a proper $\mathcal{X}$-sequence (in any order), we have the following:

 \[
val_{\mathcal{X}}(\mathcal{D})\leq val(\mathcal{D},\mathcal{A})=\left|\mathcal{D}\cup(\cup_{C\in\mathcal{D}}\mathcal{C}_{\mathcal{C}_{B},C})\right|-\left|\mathcal{D}\right|=\left|\cup_{C\in\mathcal{D}}\mathcal{C}_{\mathcal{C}_{B},C}\right|-\left|\mathcal{D}\right|\leq
\]
\[
\left|\left\{ C_{B,e}\mid e\in C\setminus B,\,C\in\mathcal{D}\right\} \right|=\left|B^{c}\right|-\left|I\right|=\rho_{\mathcal{M}^{*}}(E)-\rho_{\mathcal{M}^{*}}(\cap_{C\in\mathcal{D}}C^{c})
\]\\

\subsection{Corollary:}
If $\delta \mathcal{M}=\uparrow (2^{C(\mathcal{M})}\setminus \mathcal{S}(\mathcal{M}))$ then $val_{\mathcal{X}}$ is its rank function. In particular, we have seen that this is the case for co-rank $3$ and uniform matroids.\\

Finally, we demonstrate that, with a slight modification in its construction, $\delta	\mathcal{M}$ is not smaller than any adjoint of $\mathcal{M}$ in the weak order of matroids. Consequently, if $\delta \mathcal{M}$qualifies as an adjoint, it must be the maximal adjoint. Moreover, applying Theorem \ref{Theorem14}, we establish that if  $val_{\mathcal{X}}$ s a rank function of a matroid, then it serves as the rank function of the modified version of $\delta \mathcal{M}$. 

\subsection{Definition:}
\label{Definition14}
The small change to the combinatorial derived matroid $\delta ' \mathcal{M}$ is defined to be the matroid on the set of circuits $C(\mathcal{M})$ and dependent sets $D^{(\rho_{\mathcal{M}^*}(E))}(2^{C(\mathcal{M})}\setminus \mathcal{S})$, with $D{(k+1)}(\mathcal{A})=D(\mathcal{A}\cup (D{(k)}(\mathcal{A}))_{\leq k})$ and $D^{(0)}(\mathcal{A})=\mathcal{A}$.\\

We can see that the only change to $\delta	\mathcal{M}$ is that in $\delta \mathcal{M}$ two dependent sets can have an independent intersection, which eventually becomes dependent in a later iteration. Let $\mathcal{D}$ be a minimal, both in size and in the number of iterations required to obtain it, that remains independent in an adjoint $\mathcal{M}^\triangle$. It must be the case that $\mathcal{D}=(\mathcal{D}_1\cup \mathcal{D}_2)\setminus \left\{ C\right\} $ with $\mathcal{D}_1, \mathcal{D}_2$ dependent in $\mathcal{M}^\triangle$. Therefore, $\mathcal{D}_1\cap \mathcal{D}_2$ is dependent in $\mathcal{M}^\triangle$ and $\delta ' \mathcal{M}\nleq \mathcal{M}^\triangle$, which results in $\delta ' \mathcal{M}$ not being smaller than any adjoint of $\mathcal{M}$. We finish with a conjecture that will follow from Conjecture 1.3 in \cite{[11]}.

\subsection{Conjecture:}
Let $\mathcal{M}$ be a matroid, if there exists a maximal adjoint $\mathcal{M}^\triangle$ then $\mathcal{M}^\triangle=\delta ' \mathcal{M}$ with rank function $val_{\mathcal{X}}$.\\

%\section{The duality-polarity conjecture}

%
\section{Declarations}
\begin{enumerate}
    \item Ethical approval - not applicable.
    \item Funding - The author is supported by Horizon Europe ERC Grant number: 101045750 / Project acronym: HodgeGeoComb.
    \item Availability of data and materials - The catalog of matroids for SageMath was used for computations in the Appendix.
\end{enumerate}

\newpage
\appendix
\section{Computations}

Using Theorem \ref{Theorem8} we can construct a simple algorithm to check whether a matroid in the set of circuits of another matroid is its adjoint. We then implement this algorithm in Python using the Sage matroid package to calculate the list of adjoints of some matroids with a small number of circuits.

\begin{algorithm}
\caption{Checking if adjoint}\label{alg:cap}
\begin{algorithmic}
\Require $M,N$ are matroids and groundset($N$)=circuits($M$)
\Ensure $N$ is an adjoint of $M$
\For{B in bases(M)} \Comment{Checking if every family of fundamental circuits of a base $B$ is a base of $N$} 
            \If{$C_B$ not in bases(N)}
                \State \Return False
            \EndIf    
        \EndFor
        \For{B in bases(M)} \Comment{Checking if every fundamental circuits of $(\mathcal{L}(M^*))^{opp}$ is a circuit of $N$} 
            \For{C in circuits(M)}
                \For{$\mathcal{A}$ in bases(N)}
                    \If{$\mathcal{C}_{\mathcal{C}_{\mathcal{B}},C}$ is contained in $\mathcal{A}$}
                        \State \Return False  
                    \EndIf  
                \EndFor
            \EndFor
        \EndFor
        \State \Return True
\end{algorithmic}
\end{algorithm}

\section*{Co-rank 3 matroids}
We start by looking at co-rank $3$ examples for which the maximal adjoint $\mathcal{S}(\mathcal{M})$ is always an adjoint.

\subsection{K4}
The graphic matroid on the complete graph on four vertices has exactly two adjoints (they were already known to be adjoints in \cite{[8]}). They are the maximal and minimal adjoints and are isomorphic to the Fano and non-Fano matroids.

\subsection{Fano dual}
The dual of the fano matroid has only one adjoint, which is both the maximal and the minimal adjoint, it is isomorphic to the Fano matroid.

\subsection{non-Fano dual}
The dual of the non-fano matroid has only one adjoint, the maximal one, isomorphic to the rank $3$ ternary Dowling geometry.

\subsection{Q6}
The unique matroid representable over a field if and only if it has at least four elements has four adjoints of rank $3$ on $11$ elements:
\begin{itemize}
    \item A minimal adjoint with $136$ bases.
    \item $2$ isomorphic adjoints with $137$ bases.
    \item A maximal adjoint with $138$ bases.
\end{itemize}

\subsection{R6}
The unique matroid representable over a field if and only if it has at least three elements has sixty four adjoints of rank $3$ on $11$ elements:
\begin{itemize}
    \item A minimal adjoint with $135$ bases.
    \item $6$ isomorphic adjoints with $136$ bases.
    \item $15$ adjoints with $137$ bases of two isomorphism types, $9$ of the first type and $6$ of the seconed type.
    \item $20$ adjoints with $138$ bases of two isomorphism types, $12$ of the first type and $8$ of the seconed type.
    \item $15$ adjoints with $139$ bases of two isomorphism types, $9$ of the first type and $6$ of the seconed type.
    \item $6$ isomorphic adjoints with $140$ bases.
    \item A maximal adjoint with $141$ bases.
\end{itemize}

\subsection{P6}
The unique matroid representable over a field if and only if it has at least five elements has sixty four adjoints of rank $3$ on $13$ elements:
\begin{itemize}
    \item A minimal adjoint with $238$ bases.
    \item $6$ isomorphic adjoints with $239$ bases.
    \item $15$ adjoints with $240$ bases of two isomorphism types, $9$ of the first type and $6$ of the seconed type.
    \item $20$ adjoints with $241$ bases of two isomorphism types, $12$ of the first type and $8$ of the seconed type.
    \item $15$ adjoints with $242$ bases of two isomorphism types, $9$ of the first type and $6$ of the seconed type.
    \item $6$ isomorphic adjoints with $243$ bases.
    \item A maximal adjoint with $244$ bases.
\end{itemize}

\subsection{The dual of the represented matroid}
The dual of the first matroid from the co-rank $4$ section represented by $A$ has two adjoints. They are the maximal and minimal adjoints and are isomorphic to the Fano and non-Fano matroids.

\section*{Co-rank 4 matroids}

\subsection{Representable matroid}
The matroid represented by 
\[
A=\left(\begin{matrix}1 & 0 & 0 & 1 & 1 & 1 & 1\\
0 & 1 & 1 & 1 & 0 & 1 & 1\\
0 & 0 & 1 & 0 & 1 & 1 & 1
\end{matrix}\right)
\]
over the field of two elements has two adjoints of rank $4$ on $14$ elements. The derived matroid of $A$ which has $304$ bases and another matroid with $318$ bases. Both are not the maximal or the minimal adjoints; the first is smaller than the second in the weak order of matroids.

\subsection{AG(3,2)}
The matroid of the affine geometry of dimension three over the field of two elements has only the minimal adjoint. It is a matroid of rank $4$ on $14$ elements with $616$ bases. 


\begin{thebibliography}{BDP90}


\bibitem {[1]} White, Neil, G-C. Rota, and Neil M. White, eds. Theory of matroids. No. 26. Cambridge University Press, 1986.

\bibitem {[2]}  A. L. C. Cheung, Adjoints of a geometry, Canad. Math. Bull. 17, No. 3 (1974), 363–365.

\bibitem {[3]}  A. Bachem, A. Wanka, Matroids without adjoints, Geom. Dedicata, 29 (1989), 311–315.

\bibitem {[4]} M. Alfter, W. Kern, A. Wanka, On adjoints and dual matroids, J. Combin. Theory Ser.
B, 50 (1990), 208–213.

\bibitem {[5]} M. Alfter, W. Hochstättler, On pseudomodular matroids and adjoints, Discrete Appl.
Math. 60 (1995), 3–11.

\bibitem {[6]} R. Freij-Hollanti, R. Jurrius, O. Kuznetsova, Combinatorial derived matroid, Electron.
J. Combin. 30 (2023), Paper No. 2.8.

\bibitem {[7]}  R. Jurrius, R. Pellikaan, The coset leader and list weight enumerator, Topics in finite
fields, 229–251, Contemp. Math., 632, Amer. Math. Soc., Providence, RI, 2015.


\bibitem {[8]}  J. Oxley, S. Wang, Dependencies among dependencies in matroids, Electron. J. Combin.
26 (2019), No. 3. 46, 12pp.

\bibitem {[9]} O. Raz, Lattices of flats for symplectic matroids, arXiv:2210.15223v3 [math.CO]

\bibitem {[10]} Andreas Blass, Bruce E. Sagan. "Mobius functions of lattices".  Adv. in Math. 127 (1997), 94-123.

\bibitem {[11]} Bill Jackson and Shin-ichi Tanigawa. Maximal matroids in weak order
posets, 2021.
\bibitem {[12]} Houshan Fu and Chunming Tang and Suijie Wang.Adjoints of Matroids, arXiv:2304.08000 [math.CO]

\end{thebibliography}
\end{document}